\documentclass[11pt]{article}
\usepackage{amsfonts}
\usepackage{graphics}
\usepackage{indentfirst}
\usepackage{color}
\usepackage{cite}
\usepackage{latexsym}
\usepackage[paper=a4paper, left=2.1cm, right=2.1cm, top=2.2cm, bottom=1.5cm, headheight=5.5pt, footskip=0.8cm, footnotesep=0.8cm, centering, includefoot]{geometry}
\usepackage{amsmath}
\allowdisplaybreaks
\usepackage{amssymb}
\usepackage[dvips]{epsfig}
\usepackage{amscd}

\newtheorem{theorem}{Theorem}[section]
\newtheorem{remark}{Remark}[section]

\newtheorem{lemma}{Lemma}[section]

\DeclareMathOperator{\divv}{div}

\makeatletter
\@addtoreset{equation}{section}
\makeatother
\makeatletter
\@addtoreset{equation}{section}
\makeatother

\title{A blow-up criterion for strong solutions to three-dimensional compressible magnetohydrodynamic equations}

\author{Xin Zhong\thanks{Institute of Applied Mathematics, AMSS, Chinese Academy of Sciences, Beijing 100190,
People's Republic of China ({\tt xzhong1014@amss.ac.cn}).
}
}
\date{ }

\begin{document}
\maketitle

\begin{abstract}
We are concerned with an initial boundary value problem for the compressible magnetohydrodynamic equations with viscosity depending on the density. It is shown that for the initial density away from vacuum, the strong solution to the problem exists globally if the gradient of velocity satisfies $\|\nabla\mathbf{u}\|_{L^{2}(0,T;L^\infty)}<\infty$. Our method relies upon the delicate energy estimates and elliptic estimates.
\end{abstract}


Math Subject Classification: 76W05; 35B65

\section{Introduction}
Let $\Omega\subset\mathbb{R}^3$ be a bounded smooth domain, we consider the magnetohydrodynamic equations for compressible
fluids in $\Omega$:
\begin{align}\label{1.1}
\begin{cases}
\partial_{t}\rho+\divv(\rho\mathbf{u})=0,\\
\partial_{t}(\rho\mathbf{u})+\divv(\rho\mathbf{u}\otimes\mathbf{u})
+\nabla(A\rho^{\gamma})=
\divv(2\nu\mathfrak{D}(\mathbf{u}))
+\nabla(\lambda\divv\mathbf{u})+(\nabla\times\mathbf{H})\times\mathbf{H},\\
\partial_{t}\mathbf{H}-\mathbf{H}\cdot\nabla\mathbf{u}
+\mathbf{u}\cdot\nabla\mathbf{H}+\mathbf{H}\divv\mathbf{u}=\beta\Delta\mathbf{H},\\
\divv \mathbf{H}=0.
\end{cases}
\end{align}
Here, the unknown functions $\rho$, $\mathbf{u}, \ \mathbf{H}$ are the density, velocity of the fluid, and the magnetic field, respectively. $A$ is a positive constant, and $\gamma>1$ is the specific heat ratio. Since the concrete value of $A$ do not play a special role in our discussion, in what follows, we assume
\begin{equation*}
A=1
\end{equation*}
throughout this paper. In addition, $\mathfrak{D}(\mathbf{u})$ denotes the deformation tensor
\begin{equation*}
\mathfrak{D}(\mathbf{u})=\frac{1}{2}(\nabla\mathbf{u}+(\nabla\mathbf{u})^{tr}).
\end{equation*}

We assume that the viscosity coefficients $\nu=\nu(\rho),\lambda=\lambda(\rho)$ are twice continuously differentiable functions and satisfy
\begin{equation}\label{1.2}
\nu\geq\underline{\nu}>0,\ \mu\triangleq\nu+\lambda\geq\underline{\mu}>0.
\end{equation}
Furthermore, suppose that there is a constant $k>-1$ such that
\begin{equation}\label{1.3}
\lambda(\rho)=k\nu(\rho).
\end{equation}
The constant $\beta>0$ is the magnetic diffusive coefficient.

We study an initial boundary value problem for \eqref{1.1} with the initial condition
\begin{equation}\label{1.4}
(\rho,\mathbf{u},\mathbf{H})(0,x)=(\rho_0,\mathbf{u}_0,\mathbf{H}_0)(x),\ \ x\in\Omega,
\end{equation}
and the boundary condition
\begin{equation}\label{1.5}
\mathbf{u}=\mathbf{0},\ \mathbf{H}=\mathbf{0},\ \text{on}\ \partial\Omega.
\end{equation}
Throughout this paper, we further assume that the initial density is away from vacuum, that is,
\begin{equation}\label{1.6}
\rho_0=\rho_0(x)\geq\underline{\rho}>0,\ x\in\Omega,
\end{equation}
for a positive constant $\underline{\rho}$.

There is huge literature on the studies about the theory of
well-posedness of solutions to the Cauchy problem and initial-boundary-value problem for the system \eqref{1.1} ) in the case that the viscosity coefficients $\nu$ and $\lambda$ are both constants, refer to \cite{HW2008,HW2010,LH2015,VK1972} and references therein.
In particular, non-vacuum small perturbations of a uniform non-vacuum constant state have been shown existing globally in time and remain smooth in \cite{K1983}. However, the global existence of the strong solutions for arbitrary initial data remain largely open for the compressible magnetohydrodynamic equations with constant viscosity coefficients, but there is new progress recently \cite{LXZ2013,LSX2016}.

In the last decades, Beale-Kato-Majda \cite{BKM1984} established
a well-known blowup criterion for the three-dimensional incompressible Euler equations that a smooth solution exists globally if
\begin{equation}
\lim_{T\rightarrow T^{*}}\int_{0}^{T}\|\nabla\times\mathbf{u}\|_{L^\infty}<\infty.
\end{equation}
Later, Huang et al. \cite{HLX20111,HX2010} extended the Beale-Kato-Majda criterion for the ideal incompressible flows to the viscous compressible barotropic flows (see also \cite{FJO2010}). Furthermore, Sun-Wang-Zhang \cite{SWZ2011} established an alternative Beale-Kato-Majda blow-up criterion for the strong solution to the compressible three-dimensional Navier-Stokes equations in terms of the density instead of gradient of velocity. Suen \cite{S2013} improved the result in \cite{SWZ2011} to the compressible magnetohydrodynamic equations, which shows that the mechanism of blow-up is independent of the magnetic field. For more information on the blow-up criteria of compressible flows, we refer to \cite{HL2013,HLW2013,HLX20112,S2013,SWZ20112,WZ2014,X1998,XY2013, XZ2012} and the references therein.

Despite the progress on various blow-up criterion, and the interesting results on global well-posedness of the classical solution to the three-dimensional compressible magnetohydrodynamic equations for initial data with small total energy but possible large oscillations and containing vacuum states \cite{LXZ2013}, the behavior near vacuum of solutions to the system \eqref{1.1} remains to be one of the challenging issues for the global well-posedness of smooth solutions.

It should be noted that all the results mentioned above on the blow-up of smooth solutions of compressible flows are for constant viscosity, i.e., both $\nu$ and $\lambda$ are constant. Therefore, a natural question arises: how we can describe the mechanism of blow up to ensure the global existence of strong (or classical) solutions to the system \eqref{1.1}? Recently, Cai-Sun \cite{CS2016} obtained a blow-up criterion for the compressible Navier-Stokes equations with variable viscosity. However, the mechanism of blow-up of strong solutions to the problem \eqref{1.1} with variable viscosity is still
unknown. In fact, this is the main aim of this paper.

Our main result reads as follows:
\begin{theorem}\label{thm1.1}
Assume that the initial data $(\rho_0,\mathbf{u}_0,\mathbf{H}_0)$ satisfy \eqref{1.6} and
\begin{align}\label{A}
(\rho_0,\mathbf{u}_0,\mathbf{H}_0)\in H^{3}(\Omega)
\end{align}
with the compatibility conditions on $\partial\Omega$,
\begin{equation}\label{A2}
\mathbf{u}_0=\mathbf{0},\ \mathbf{H}_0=\mathbf{0},\
\nabla\rho_{0}^{\gamma}=\divv(2\nu(\rho_0)\mathfrak{D}(\mathbf{u_0}))
+\nabla(\lambda(\rho_0)\divv\mathbf{u}_0).
\end{equation}
Let $(\rho,\mathbf{u},\mathbf{H})\in C([0,T^{*}),H^3(\Omega))$ be a strong solution to the problem \eqref{1.1} with the conditions \eqref{1.2}--\eqref{1.5}.
If $T^{*}<\infty$ is the maximal time of existence for that solution,
then we have
\begin{align}\label{B}
\lim_{T\rightarrow T^{*}}\|\nabla\mathbf{u}\|_{L^{2}(0,T;L^\infty)}=\infty.
\end{align}
\end{theorem}

Several remarks are in order.
\begin{remark}\label{re1.1}
The local existence of a strong solution with initial data as in Theorem \ref{thm1.1} can be established in the same manner as \cite{V1982}. Hence the maximal time $T^{*}$ is well-defined.
\end{remark}
\begin{remark}\label{re1.2}
It is worth noting that the blow-up criteria \eqref{B} is independent of the magnetic field, and is just the same as that of the compressible Navier-Stokes equations \cite[Theorem 1.4]{CS2016}. So Theorem \ref{thm1.1} may be viewed as an extension of \cite[Theorem 1.4]{CS2016}.
\end{remark}
\begin{remark}\label{re1.3}
The techniques in this paper can be easily adapted to the periodic case. And the same criterion will be derived.
\end{remark}
\begin{remark}\label{re1.4}
Similar ideas can be applied to study the case on the viscous, compressible, and heat conducting magnetohydrodynamic equations with variable viscosity.
This will be reported in a forthcoming paper \cite{Z2016}.
\end{remark}

We now make some comments on the analysis of this paper. We mainly make use of continuation argument to show Theorem \ref{thm1.1}. Roughly speaking, under the assumption that \eqref{B} were false, one needs estimate $L^{\infty}(0,T;H^3)$-norm for any $T\in[0,T^{*})$ on strong solutions $(\rho,\mathbf{u},\mathbf{H})$ to the system \eqref{1.1}. It should be pointed out that, on one hand, the crucial techniques of proofs in \cite{S2013,XZ2012} cannot be adapted to the situation treated here, since their arguments only hold true for the case of absence of boundary. On the other hand, due to the presence of variable coefficients, it seems nontrivial for derivation of higher order energy estimates of solutions. To this end, we try to adapt some basic ideas used in \cite{CS2016,SWZ2011}, where they investigated the mechanism of blow-up of strong solutions to the compressible Navier-Stokes(-Fourier) equations. However, compared with \cite{CS2016,SWZ2011}, for the compressible MHD equations treated here, the nonlinearities of the magnetic field itself and its interaction with  the velocity field of the fluid, such as $\mathbf{u}\cdot\nabla\mathbf{H}$, will bring out some new difficulties.

To overcome these difficulties stated  above, some new elaborate estimates are needed. First, we try to obtain higher integrability of the magnetic field (see Lemma \ref{lem33}), inspired by \cite{HX2005} (see also \cite{HL2013,XZ2012}), we multiply \eqref{1.1}$_3$ by $q|\mathbf{H}|^{q-2}\mathbf{H}$ and
thus get $L^\infty(0,T;L^q)$-norm of the magnetic field
$\mathbf{H}$, which is crucial in deriving the higher order estimates of the solutions. Next, we attempt to show the estimates on the $L^\infty(0,T;L^2)$-norm of the gradients of velocity and magnetic. Motivated by \cite{HL2013}, multiplying \eqref{1.1}$_2$ by the time derivatives $\mathbf{u}_t$, the key point is to control the term $\|\nabla^2\mathbf{u}\|_{L^2}$ (see \eqref{3.11} and \eqref{3.14}). Motivated by \cite{SWZ2011}, using some facts on Lam{\'e} system with constant coefficient (see Lemma \ref{lem23}), we succeed in bounding $\|\nabla^2\mathbf{u}\|_{L^2}$ in terms by $\|\nabla^2\mathbf{H}\|_{L^2}$ (see \eqref{b}). Then, multiplying \eqref{1.1}$_3$ by $\Delta\mathbf{H}$, the term $\|\nabla^2\mathbf{H}\|_{L^2}$ can be absorbed (see \eqref{3.20}). On the other hand, the usual $L^{\infty}(0,T;L^2)$-norm of $\mathbf{u}_t$ and $\mathbf{H}_t$ cannot be directly estimated as usual by taking inner product \eqref{1.1}$_2$ and \eqref{1.1}$_3$ with $\Delta\mathbf{u}$ and $\Delta\mathbf{H}$ respectively due to the presence of boundary. Motivated by \cite{CS2016}, differentiating \eqref{1.1}$_2$ and \eqref{1.1}$_3$ with respect to $t$, and multiplying the resulting equations by $\mathbf{u}_t$ and $\mathbf{H}_t$, respectively, (see \eqref{3.23} and \eqref{4.2}) and thus obtain $L^\infty(0,T;L^2)$-norm of $\sqrt{\rho}\mathbf{u}_t$ and $\mathbf{H}_t$ (see \eqref{4.6}). This together with some analysis on the magnetic field $\mathbf{H}$ and Lam{\'e} system indicate the $H^3$-bound on spatial one of the velocity (see \eqref{z4})
can be estimated in terms of $\|\nabla^2\rho\|_{L^2}$,
which in turn implies the bound on the $L^\infty(0,T;H^2)$-norm of the the density (see \eqref{z2}). With the a priori estimates stated above at hand, the desired $L^{\infty}(0,T;H^2)$-bound of the velocity and the magnetic field can be deduced (see \eqref{4.12} and \eqref{4.13}).
Finally, $L^{\infty}(0,T;H^3)$-norm on $(\rho,\mathbf{u},\mathbf{H})$ (see \eqref{l5}, \eqref{l7} and \eqref{l8}) is obtained, such a derivation is similar as the estimate on the $L^{\infty}(0,T;H^2)$-norm on $(\rho,\mathbf{u},\mathbf{H})$, see Lemma \ref{lem36} and its proof.

The rest of this paper is organized as follows. In Section \ref{sec2}, we collect some elementary facts and inequalities that will be used later. Section \ref{sec3} is devoted to the proof of Theorem \ref{thm1.1}.

\section{Preliminaries}\label{sec2}

In this section, we will recall some known facts and elementary inequalities that will be used frequently later.

We begin with the following Gronwall's inequality, which plays a central role in proving a priori estimates on strong solutions $(\rho,\mathbf{u},\mathbf{H})$.
\begin{lemma}\label{lem21}
Suppose that $h$ and $r$ are integrable on $(a, b)$ and nonnegative a.e. in $(a, b)$. Further assume that $y\in C[a, b], y'\in L^1(a, b)$, and
\begin{equation*}
y'(t)\leq h(t)+r(t)y(t)\ \ \text{for}\ a.e\ t\in(a,b).
\end{equation*}
Then
\begin{equation*}
y(t)\leq \left[y(a)+\int_{a}^{t}h(s)\exp\left(-\int_{a}^{s}r(\tau)d\tau\right)ds\right]
\exp\left(\int_{a}^{t}r(s)ds\right),\ \ t\in[a,b].
\end{equation*}
\end{lemma}
{\it Proof.}
See  \cite[pp. 12--13]{T2006}.  \hfill $\Box$

Next, the following well-known inequalities will be frequently used later.
\begin{lemma}\label{lem22}
For $p\in[2, 6],q,m\in[1,\infty),\alpha\in(0,m)$, $\varepsilon>0,a,b\in\mathbb{R}$, and $\theta\in(0,1)$, it holds that
\begin{equation*}
\begin{split}
& \|f\|_{L^m}\leq \|f\|_{L^{\alpha q}}^{\frac{\alpha}{m}}
\|f\|_{L^{(m-\alpha)q}}^{1-\frac{\alpha}{m}},\ \ \ \ \ \ \ \ \ \ \ \ \ \ \text{(H{\"o}lder's inequality)} \\
& |ab|\leq \varepsilon|a|^{\frac{1}{\theta}}
+\left(\frac{\theta}{\varepsilon}\right)^{\frac{\theta}{1-\theta}}
(1-\theta)|b|^{\frac{1}{1-\theta}},\ \ \ \ \text{(Young's inequality)}
\end{split}
\end{equation*}
and
\begin{equation*}
\|g\|_{L^p}\leq C(p,\Omega)\|g\|_{H^1}\ \ \text{for}\ \ g\in H^{1}(\Omega).
\ \ \text{(Sobolev's inequality)}
\end{equation*}
\end{lemma}
{\it Proof.}
See  \cite[Chapter 2]{LU1968}.  \hfill $\Box$

Finally, we give some regularity results for the following Lam{\'e} system with Dirichlet boundary condition
\begin{align}\label{2.1}
\begin{cases}
\Delta\mathbf{U}+\nabla\divv\mathbf{U}=\mathbf{F},\ \ x\in\Omega,\\
\mathbf{U}=\mathbf{0},\ \ x\in\partial\Omega.
\end{cases}
\end{align}
\begin{lemma}\label{lem23}
Let $q\in(1,\infty)$ and U be a solution of \eqref{2.1}. There exists a constant $C$ depending only on $q$
and $\Omega$ such that the following estimates hold:
\begin{itemize}
\item [$\bullet$] If $\mathbf{F}\in L^{q}(\Omega)$, then
\begin{equation*}
\|\mathbf{U}\|_{W^{2,q}}\leq C\|\mathbf{F}\|_{L^q};
\end{equation*}
\item [$\bullet$] If $\mathbf{F}\in W^{-1,q}(\Omega)$ (i.e., $\mathbf{F}= \divv\mathbf{f}$ with $\mathbf{f}=(f_{ij})_{3\times3}, f_{ij}\in L^q(\Omega)$), then
\begin{equation*}
\|\mathbf{U}\|_{W^{1,q}}\leq C\|\mathbf{f}\|_{L^q};
\end{equation*}
\item [$\bullet$] If $\mathbf{F}= \divv\mathbf{f}$ with $f_{ij}=\partial_{k}h_{ij}^{k}$ and $h_{ij}^{k}\in W_{0}^{1,q}(\Omega)$) for $i,j,k=1,2,3$, then
\begin{equation*}
\|\mathbf{U}\|_{L^{q}}\leq C\|\mathbf{h}\|_{L^q}.
\end{equation*}
\end{itemize}
\end{lemma}
{\it Proof.}
See  \cite[Proposition 2.1]{SWZ2011}.  \hfill $\Box$

\section{Proof of Theorem \ref{thm1.1}}\label{sec3}

Let $(\rho,\mathbf{u},\mathbf{H})$ be a strong solution described in Theorem \ref{thm1.1}. Suppose that \eqref{B} were false, that is, there exists a constant $M_0>0$ such that
\begin{equation}\label{3.1}
\lim_{T\rightarrow T^*}\|\nabla\mathbf{u}\|_{L^{2}(0,T;L^\infty)}\leq M_0<\infty.
\end{equation}
For simplicity, in what follows, we denote by
\begin{equation*}
\int\cdot\text{d}x=\int_{\Omega}\cdot\text{d}x,
\end{equation*}
and $C,C_1,C_2$ stand for generic positive constants depending only on $M_0,\beta,T^{*}$, and the initial data.
\begin{lemma}\label{lem31}
There exist two positive constants $C_1,C_2$ depending only on $M_0$ and the upper and lower bounds of $\rho_0$ such that
\begin{equation}\label{3.2}
C_2\leq\rho(t,x)\leq C_1\ \ \text{for all}\ \ (t,x)\in[0,T]\times\Omega.
\end{equation}
\end{lemma}
{\it Proof.}
In accordance with the continuity equation \eqref{1.1}$_1$, the density $\rho$ computed along characteristics satisfies
\begin{equation*}
\frac{d}{dt}\rho(t,\mathbf{X}(t))
=-\rho(t,\mathbf{X}(t))\divv\mathbf{u}(t,\mathbf{X}(t)),\
\mathbf{X}(0)=x,\ t\in[0,T].
\end{equation*}
Thus, a direct application of Gronwall's inequality yields that for any $(t,x)\in[0,T]\times\Omega$,
\begin{align*}
\inf_{x\in\Omega}\rho_0(x)
\exp\left(-\int_{0}^{T}\|\divv\mathbf{u}\|_{L^\infty}dt\right)
\leq\rho(t,x)\leq\sup_{x\in\Omega}\rho_0(x)
\exp\left(\int_{0}^{T}\|\divv\mathbf{u}\|_{L^\infty}dt\right),
\end{align*}
which combined with \eqref{3.1} leads to the desired result \eqref{3.2}.
\hfill $\Box$
\begin{remark}\label{re3.1}
Since $\nu,\lambda$ are both twice continuously differentiable functions, we deduce from \eqref{3.2} that
\begin{equation}\label{z}
\|(\nu',\nu'',\lambda',\lambda'')\|_{L^{\infty}(0,T;L^\infty(\Omega))}<\infty.
\end{equation}
It will be used frequently in our later analysis.
\end{remark}

\begin{lemma}\label{lem32}
There exists a positive constant $C$ depending only on $\|\sqrt{\rho_0}\mathbf{u}_0\|_{L^2}$, $\|\rho_0\|_{L^\gamma}$, and $\|\mathbf{H}_0\|_{L^2}$ such that for any $T\in[0,T^*)$, it holds that
\begin{equation}\label{3.3}
\sup_{0\leq t\leq T}\int\left(\rho|\mathbf{u}|^2+\rho^{\gamma}+|\mathbf{H}|^2\right)dx
+\int_{0}^{T}\int\left(|\nabla\mathbf{u}|^2
+|\nabla\mathbf{H}|^2\right)dxdt
\leq C.
\end{equation}
\end{lemma}
{\it Proof.}
Multiplying \eqref{1.1}$_2$ and \eqref{1.1}$_3$ by $\mathbf{u}$ and $\mathbf{H}$ respectively, then adding the two resulting equations together, and integrating over $\Omega$, we can derive that
\begin{equation}\label{3.4}
\frac{d}{dt}\int\left(\frac{1}{2}\rho|\mathbf{u}|^2
+\frac{\rho^\gamma}{\gamma-1}+\frac{1}{2}|\mathbf{H}|^2\right)dx
+\int\left(2\nu|\mathfrak{D}(\mathbf{u})|^2+\lambda|\divv\mathbf{u}|^2
+\beta|\nabla\mathbf{H}|^2\right)dx=0,
\end{equation}
where integration by parts was applied. It implies that the inequality \eqref{3.3} holds and consequently completes the proof.
\hfill $\Box$

Inspired by \cite{HX2005},
we have the following higher integrability of the magnetic field $\mathbf{H}$.
\begin{lemma}\label{lem33}
Under the condition \eqref{3.1}, it holds that for any $q\in[2,\infty)$,
\begin{equation}\label{3.5}
\sup_{0\leq t\leq T}\int|\mathbf{H}|^{q}dx
+\int_{0}^{T}\int|\mathbf{H}|^{q-2}|\nabla\mathbf{H}|^2dxdt\leq C.
\end{equation}
\end{lemma}
{\it Proof.}
Multiplying \eqref{1.1}$_3$ by $q|\mathbf{H}|^{q-2}\mathbf{H}$ and integrating the resulting equations over $\Omega$, we obtain
\begin{align}\label{3.6}
& \frac{d}{dt}\int|\mathbf{H}|^qdx
+\beta\int\left(q|\mathbf{H}|^{q-2}|\nabla\mathbf{H}|^2
+q(q-2)|\mathbf{H}|^{q-2}|\nabla|\mathbf{H}||^2\right)dx
\nonumber\\
& = -\int q|\mathbf{H}|^{q-2}\left(\mathbf{H}\cdot\nabla\mathbf{H}\cdot\mathbf{u}
-\frac{q-1}{2}\mathbf{u}\cdot\nabla|\mathbf{H}|^2\right)dx
\nonumber\\
& \quad
-\frac{q(q-2)}{2}\int|\mathbf{H}|^{q-4}(\mathbf{H}\cdot\nabla|\mathbf{H}|^2)
(\mathbf{u}\cdot\mathbf{H})dx
\nonumber\\
& \leq \frac{\beta}{2}\int q|\mathbf{H}|^{q-2}|\nabla\mathbf{H}|^2dx
+C(q)\int|\mathbf{u}|^2|\mathbf{H}|^qdx
\nonumber\\
& \leq \frac{\beta}{2}\int q|\mathbf{H}|^{q-2}|\nabla\mathbf{H}|^2dx
+C\|\mathbf{u}\|_{L^\infty}^2\int|\mathbf{H}|^qdx
\nonumber\\
& \leq \frac{\beta}{2}\int q|\mathbf{H}|^{q-2}|\nabla\mathbf{H}|^2dx
+C\|\nabla\mathbf{u}\|_{L^\infty}^2\int|\mathbf{H}|^qdx,
\end{align}
where in the second inequality one has used
\begin{equation}\label{a}
\|\mathbf{u}\|_{L^\infty}
\leq C
\|\nabla\mathbf{u}\|_{L^p}
\leq C
\|\nabla\mathbf{u}\|_{L^\infty}\ \ \text{for any}\ \ p>3
\end{equation}
according to $\mathbf{u}|_{\partial\Omega}=\mathbf{0}$ and Sobolev's embedding theorem. Then we derive \eqref{3.5} directly after using Gronwall's inequality and \eqref{3.1}. Thus we finish the proof of Lemma \ref{lem33}. \hfill $\Box$

The following lemma gives the estimates on the spatial gradients of both the velocity and the magnetic field, which are crucial for deriving the higher order estimates of the solution.
\begin{lemma}\label{lem34}
Under the condition \eqref{3.1}, it holds that for any $T\in[0,T^*)$,
\begin{equation}\label{3.7}
\sup_{0\leq t\leq T}\int\left(|\nabla\mathbf{u}|^{2}+|\nabla\rho|^{2}+|\nabla\mathbf{H}|^{2}\right)dx
+\int_{0}^{T}\int\left(\rho|\partial_t\mathbf{u}|^{2}
+|\nabla^2\mathbf{H}|^2\right)dxdt \leq C.
\end{equation}
\end{lemma}
{\it Proof.}
Note that the momentum equation \eqref{1.1}$_2$ can be rewritten as
\begin{equation}\label{3.8}
\rho\partial_{t}\mathbf{u}-\nu\Delta\mathbf{u}-\mu\nabla\divv\mathbf{u}
=\mathbf{f},
\end{equation}
where
\begin{equation}\label{6.1}
\mathbf{f}=\nabla\nu\cdot\nabla\mathbf{u}+\nabla\mathbf{u}\cdot\nabla\nu
+\divv\mathbf{u}\nabla\lambda-\rho\mathbf{u}\cdot\nabla\mathbf{u}
-\nabla\rho^\gamma+(\nabla\times\mathbf{H})\times\mathbf{H}.
\end{equation}
By Lemma \ref{lem22} and \eqref{3.5}, we have
\begin{align}\label{6.2}
\|(\nabla\times\mathbf{H})\times\mathbf{H}\|_{L^2}
& \leq C\||\mathbf{H}||\nabla\mathbf{H}|\|_{L^2} \nonumber \\
& \leq C\|\mathbf{H}\|_{L^6}\|\nabla\mathbf{H}\|_{L^2}^{\frac{1}{2}}
\|\nabla\mathbf{H}\|_{L^6}^{\frac{1}{2}} \nonumber
\\
& \leq C\|\nabla\mathbf{H}\|_{L^2}^{\frac{1}{2}}
\|\nabla\mathbf{H}\|_{H^1}^{\frac{1}{2}}.
\end{align}
Then it follows from \eqref{6.1}, \eqref{6.2}, \eqref{3.2}, \eqref{z}, and \eqref{a} that
\begin{align}\label{3.9}
\|\mathbf{f}\|_{L^2}
& \leq C\||\nabla\mathbf{u}||\nu'||\nabla\rho|\|_{L^2}
+\||\nabla\mathbf{u}||\lambda'||\nabla\rho|\|_{L^2}
+\||\rho||\mathbf{u}||\nabla\mathbf{u}|\|_{L^2}
\nonumber \\
& \quad
+C\|\rho^{\gamma-1}|\nabla\rho|\|_{L^2}
+\|(\nabla\times\mathbf{H})\times\mathbf{H}\|_{L^2}
\nonumber \\
& \leq C\|(\nu',\lambda')\|_{L^\infty}
\|\nabla\mathbf{u}\|_{L^\infty}\|\nabla\rho\|_{L^2}
+\|\rho\|_{L^\infty}\|\mathbf{u}\|_{L^\infty}\|\nabla\mathbf{u}\|_{L^2}
\nonumber \\
& \quad +C\|\rho\|_{L^\infty}^{\gamma-1}\|\nabla\rho\|_{L^2}
+C\|\nabla\mathbf{H}\|_{L^2}^{\frac{1}{2}}
\|\nabla\mathbf{H}\|_{H^1}^{\frac{1}{2}}
\nonumber \\
& \leq C\left(1+\|\nabla\mathbf{u}\|_{L^\infty}\right)
\left(\|\nabla\mathbf{u}\|_{L^2}+\|\nabla\rho\|_{L^2}\right)
+C\|\nabla\mathbf{H}\|_{L^2}^{\frac{1}{2}}
\|\nabla\mathbf{H}\|_{H^1}^{\frac{1}{2}}.
\end{align}

Multiplying \eqref{3.8} by $\partial_{t}\mathbf{u}$ and integrating the resulting equations over $\Omega$ give rise to
\begin{align}\label{3.10}
& \int\rho|\partial_{t}\mathbf{u}|^2dx
+\frac{1}{2}\frac{d}{dt}\int\left(\nu|\nabla\mathbf{u}|^2+\mu|\divv\mathbf{u}|^2
\right)dx\nonumber\\
& = \int\mathbf{f}\cdot\partial_{t}\mathbf{u}dx
+\frac{1}{2}\int\partial_{t}\nu|\nabla\mathbf{u}|^2dx
-\int\nabla\nu\cdot\nabla\mathbf{u}\cdot\partial_{t}\mathbf{u}dx
\nonumber\\
& \quad
+\frac{1}{2}\int\partial_{t}\mu|\divv\mathbf{u}|^2dx
-\int\partial_{t}\mathbf{u}\nabla\mu\divv\mathbf{u}dx
\triangleq \sum_{j=1}^{5}I_j.
\end{align}
Employing straightforward calculations, we can estimate each term $I_j$ as follows
\begin{align*}
|I_{1}|\leq & \|\rho^{-\frac{1}{2}}\|_{L^\infty}
\|\mathbf{f}\|_{L^2}\|\sqrt{\rho}\partial_{t}\mathbf{u}\|_{L^2}
\leq \frac{1}{6}\|\sqrt{\rho}\partial_{t}\mathbf{u}\|_{L^2}^2
+C\|\mathbf{f}\|_{L^2}^2 \\
\leq & \frac{1}{6}\|\sqrt{\rho}\partial_{t}\mathbf{u}\|_{L^2}^2
+C\left(1+\|\nabla\mathbf{u}\|_{L^\infty}^2\right)
\left(1+\|\nabla\mathbf{u}\|_{L^2}^2+\|\nabla\rho\|_{L^2}^2\right)
+C\|\nabla\mathbf{H}\|_{L^2}\|\nabla\mathbf{H}\|_{H^1};\\
|I_{2}|\leq & \frac{1}{2}\|\nu'\|_{L^\infty}\|\partial_{t}\rho\|_{L^2}\|\nabla\mathbf{u}\|_{L^2} \|\nabla\mathbf{u}\|_{L^\infty}
\leq C\|\partial_{t}\rho\|_{L^2}\|\nabla\mathbf{u}\|_{L^2} \|\nabla\mathbf{u}\|_{L^\infty}\\
\leq & C\|\nabla\mathbf{u}\|_{L^\infty}^{2}
\left(1+\|\nabla\mathbf{u}\|_{L^2}^2+\|\nabla\rho\|_{L^2}^2\right);\\
|I_{3}|\leq & \|\rho^{-\frac{1}{2}}\|_{L^\infty} \|\nu'\|_{L^\infty}\|\nabla\rho\|_{L^2}
\|\sqrt{\rho}\partial_{t}\mathbf{u}\|_{L^2} \|\nabla\mathbf{u}\|_{L^\infty}
\leq\frac{1}{6}\|\sqrt{\rho}\partial_{t}\mathbf{u}\|_{L^2}^2
+C\|\nabla\mathbf{u}\|_{L^\infty}^2\|\nabla\rho\|_{L^2}^2;\\
|I_{4}|\leq & \frac{1}{2}\|\mu'\|_{L^\infty}\|\partial_{t}\rho\|_{L^2}\|\nabla\mathbf{u}\|_{L^2} \|\nabla\mathbf{u}\|_{L^\infty}
\leq C\|\partial_{t}\rho\|_{L^2}\|\nabla\mathbf{u}\|_{L^2} \|\nabla\mathbf{u}\|_{L^\infty}\\
\leq & C\|\nabla\mathbf{u}\|_{L^\infty}^{2}
\left(1+\|\nabla\mathbf{u}\|_{L^2}^2+\|\nabla\rho\|_{L^2}^2\right);\\
|I_{5}|\leq & \|\rho^{-\frac{1}{2}}\|_{L^\infty} \|\mu'\|_{L^\infty}\|\nabla\rho\|_{L^2}
\|\sqrt{\rho}\partial_{t}\mathbf{u}\|_{L^2} \|\nabla\mathbf{u}\|_{L^\infty}
\leq\frac{1}{6}\|\sqrt{\rho}\partial_{t}\mathbf{u}\|_{L^2}^2
+C\|\nabla\mathbf{u}\|_{L^\infty}^2\|\nabla\rho\|_{L^2}^2,
\end{align*}
where we have used
\begin{equation}\label{x}
\partial_{t}\rho=-\mathbf{u}\cdot\nabla\rho-\rho\divv\mathbf{u},\
\|\partial_{t}\rho\|_{L^2}
\leq C\|\nabla\mathbf{u}\|_{L^\infty}(1+\|\nabla\rho\|_{L^2})\quad (\text{by} \ \eqref{a}\ \text{and}\ \eqref{3.2})
\end{equation}
in showing the estimates of $I_2$ and $I_4$, and \eqref{3.2} in obtaining the estimates of $I_1,I_3,$ and $I_5$.
Substituting the above estimates on $I_j$ into \eqref{3.10}, we infer from \eqref{3.5} that
\begin{align}\label{3.11}
& \frac{d}{dt}\int\left(\nu|\nabla\mathbf{u}|^2+\mu|\divv\mathbf{u}|^2
\right)dx
+\int\rho|\partial_{t}\mathbf{u}|^2dx \nonumber\\
& \leq C\left(1+\|\nabla\mathbf{u}\|_{L^\infty}^2\right)
\left(1+\|\nabla\mathbf{u}\|_{L^2}^2+\|\nabla\rho\|_{L^2}^2\right)
+C\|\nabla\mathbf{H}\|_{L^2}\|\nabla\mathbf{H}\|_{H^1}
\nonumber\\
& \leq C\left(1+\|\nabla\mathbf{u}\|_{L^\infty}^2\right)
\left(1+\|\nabla\mathbf{u}\|_{L^2}^2+\|\nabla\rho\|_{L^2}^2\right)
+C\|\nabla\mathbf{H}\|_{L^2}^2+C\|\nabla^2\mathbf{H}\|_{L^2}^2.
\end{align}

Now we shall need to estimate $\|\nabla\rho\|_{L^2}$. To this end, taking spatial derivative $\nabla$ on the transport equation \eqref{1.1}$_1$ leads to
\begin{equation*}
\partial_{t}\nabla\rho+\mathbf{u}\cdot\nabla^2\rho+\nabla\mathbf{u}\cdot\nabla\rho
+\divv\mathbf{u}\nabla\rho+\rho\nabla\divv\mathbf{u}=\mathbf{0}.
\end{equation*}
Then standard energy methods yields
\begin{align}\label{3.14}
\frac{d}{dt}\int|\nabla\rho|^2dx
\leq C\left(1+\|\nabla\mathbf{u}\|_{L^\infty}^2\right)\|\nabla\rho\|_{L^2}^2
+C\|\nabla^2\mathbf{u}\|_{L^2}^2.
\end{align}
To control the last term in \eqref{3.14}, we see that, by \eqref{1.2} and \eqref{1.3}, \eqref{1.1}$_2$ can be rewritten as the following Lam{\'e} system
\begin{align}\label{x2}
\begin{cases}
\Delta\mathbf{u}+\left(k+1\right)\nabla\divv\mathbf{u}
=\frac{\mathbf{F}}{\nu},\ \ x\in\Omega,\\
\mathbf{u}=\mathbf{0},\ \ x\in\partial\Omega,
\end{cases}
\end{align}
with
\begin{equation}\label{x3}
\mathbf{F}=\rho\partial_{t}\mathbf{u}
+\rho\mathbf{u}\cdot\nabla\mathbf{u}
+\nabla\rho^\gamma
-\nabla\nu\cdot\nabla\mathbf{u}
-\nabla\mathbf{u}\cdot\nabla\nu
-\divv\mathbf{u}\nabla\lambda
-(\nabla\times\mathbf{H})\times\mathbf{H}.
\end{equation}
Then we obtain from Lemma \ref{lem23}, H{\"o}lder's inequality, and \eqref{3.9} that
\begin{align}\label{b}
\|\nabla^2\mathbf{u}\|_{L^2}
& \leq C\|\mathbf{F}\|_{L^2}
\nonumber\\
& \leq \frac{1}{6}\|\sqrt{\rho}\partial_{t}\mathbf{u}\|_{L^2}+
C\left(1+\|\nabla\mathbf{u}\|_{L^\infty}\right)
\left(1+\|\nabla\mathbf{u}\|_{L^2}+\|\nabla\rho\|_{L^2}\right)
+C\|\nabla\mathbf{H}\|_{L^2}^{\frac{1}{2}}\|\nabla\mathbf{H}\|_{H^1}^{\frac{1}{2}}.
\end{align}
Inserting \eqref{b} into \eqref{3.14}, we get
\begin{align}\label{c}
\frac{d}{dt}\int|\nabla\rho|^2dx
& \leq \frac{1}{4}\|\sqrt{\rho}\partial_{t}\mathbf{u}\|_{L^2}^2+
C\left(1+\|\nabla\mathbf{u}\|_{L^\infty}^2\right)
\left(1+\|\nabla\mathbf{u}\|_{L^2}^2+\|\nabla\rho\|_{L^2}^2\right)
\nonumber\\
& \quad
+C\|\nabla\mathbf{H}\|_{L^2}\|\nabla\mathbf{H}\|_{H^1} \nonumber\\
& \leq \frac{1}{4}\|\sqrt{\rho}\partial_{t}\mathbf{u}\|_{L^2}^2+
C\left(1+\|\nabla\mathbf{u}\|_{L^\infty}^2\right)
\left(1+\|\nabla\mathbf{u}\|_{L^2}^2+\|\nabla\rho\|_{L^2}^2\right)
\nonumber\\
& \quad
+C\|\nabla\mathbf{H}\|_{L^2}^2+C\|\nabla^2\mathbf{H}\|_{L^2}^2.
\end{align}
Substituting \eqref{c} into \eqref{3.11}, we obtain
\begin{align}\label{3.18}
& \frac{d}{dt}\int\left(\nu|\nabla\mathbf{u}|^2
+\mu|\divv\mathbf{u}|^2\right)dx
+ \frac{d}{dt}\int|\nabla\rho|^2dx
+\frac{1}{2}\int\rho|\partial_{t}\mathbf{u}|^2dx
\nonumber\\
& \leq C\left(1+\|\nabla\mathbf{u}\|_{L^\infty}^2\right)
\left(1+\|\nabla\mathbf{u}\|_{L^2}^2+\|\nabla\rho\|_{L^2}^2
+\|\nabla\mathbf{H}\|_{L^2}^2\right)
+C\|\nabla^2\mathbf{H}\|_{L^2}^2.
\end{align}

Next, multiplying \eqref{1.1}$_3$ by $\Delta\mathbf{H}$, and integrating by parts over $\Omega$, we have
\begin{equation*}
\begin{split}
\frac{d}{dt}\int|\nabla\mathbf{H}|^2dx
+2\beta\int|\Delta\mathbf{H}|^2dx
& \leq C\int|\mathbf{u}||\nabla\mathbf{H}||\Delta\mathbf{H}|dx
+C\int|\mathbf{H}||\nabla\mathbf{u}||\Delta\mathbf{H}|dx
\nonumber\\
& \leq \beta\|\Delta\mathbf{H}\|_{L^2}^2
+C\||\mathbf{u}||\nabla\mathbf{H}|\|_{L^2}^2
+C\||\mathbf{H}||\nabla\mathbf{u}|\|_{L^2}^2,
\end{split}
\end{equation*}
which combined with \eqref{3.5} and \eqref{a} yields
\begin{align}\label{3.19}
\frac{d}{dt}\int|\nabla\mathbf{H}|^2dx
+\beta\int|\Delta\mathbf{H}|^2dx
& \leq C\||\mathbf{u}||\nabla\mathbf{H}|\|_{L^2}^2
+C\||\mathbf{H}||\nabla\mathbf{u}|\|_{L^2}^2 \nonumber\\
& \leq C\|\mathbf{u}\|_{L^\infty}^2\|\nabla\mathbf{H}\|_{L^2}^2
+C\|\mathbf{H}\|_{L^2}^2\|\nabla\mathbf{u}\|_{L^\infty}^2
\nonumber\\
& \leq C\|\nabla\mathbf{u}\|_{L^\infty}^2\|\nabla\mathbf{H}\|_{L^2}^2
+C\|\nabla\mathbf{u}\|_{L^\infty}^2.
\end{align}
Adding \eqref{3.19} multiplied by $C_3$ large enough to \eqref{3.18} leads to
\begin{align}\label{3.20}
& A'(t)+\frac{1}{2}\int\rho|\partial_{t}\mathbf{u}|^2dx
+\widetilde{C}\int|\nabla^2\mathbf{H}|^2dx
\nonumber\\
& \leq C\left(1+\|\nabla\mathbf{u}\|_{L^\infty}^2\right)
\left(1+\|\nabla\mathbf{u}\|_{L^2}^2+\|\nabla\rho\|_{L^2}^2
+\|\nabla\mathbf{H}\|_{L^2}^2\right)
\end{align}
for some $\widetilde{C}>0$, where, due to \eqref{1.2},
\begin{equation*}
\begin{split}
A(t) & \triangleq \int\left(\nu|\nabla\mathbf{u}|^2
+\mu|\divv\mathbf{u}|^2\right)dx
+ \int|\nabla\rho|^2dx
+ C_3\int|\nabla\mathbf{H}|^2dx+1 \\
& \geq C_{4}\int|\nabla\mathbf{u}|^2dx
+ \int|\nabla\rho|^2dx
+ C_3\int|\nabla\mathbf{H}|^2dx+1
\\
& \geq \overline{C}\left(1+\|\nabla\mathbf{u}\|_{L^2}^2+\|\nabla\rho\|_{L^2}^2
+\|\nabla\mathbf{H}\|_{L^2}^2\right)
\end{split}
\end{equation*}
for some $C_4>0$ and $\overline{C}\triangleq\min\{1,C_3,C_4\}$.
Then, the desired \eqref{3.7} is a direct consequence of \eqref{3.20}, \eqref{3.1}, and Gronwall's inequality. This completes the proof of Lemma \ref{lem34}.
\hfill $\Box$
\begin{remark}\label{re3.1}
It follows from \eqref{b}, \eqref{3.7}, and \eqref{3.1}, we infer that
\begin{equation}\label{8.1}
\int_{0}^{T}\int|\nabla^2\mathbf{u}|^2dxdt \leq C.
\end{equation}
\end{remark}

In the following two lemmas, we will focus on the estimates of $L^{\infty}(0,T;H^2)$-norm and $L^{\infty}(0,T;H^3)$-norm of strong solutions $(\rho,\mathbf{u},\mathbf{H})$, which are needed to guarantee the extension of the maximal time $T^{*}$.
\begin{lemma}\label{lem35}
Under the condition \eqref{3.1}, it holds that for any $T\in[0,T^*)$,
\begin{equation}\label{123}
\sup_{0\leq t\leq T}\left(\|\mathbf{u}\|_{H^2}
+\|\rho\|_{H^2}+\|\mathbf{H}\|_{H^2}\right) \leq C.
\end{equation}
\end{lemma}
{\it Proof.}
Differentiating \eqref{1.1}$_2$ with respect to $t$, we obtain that
\begin{align}\label{3.21}
& \rho\mathbf{u}_{tt}+\rho\mathbf{u}\cdot\nabla\mathbf{u}_{t}
-\divv\left(2\nu\mathfrak{D}(\mathbf{u}_t)
+\lambda\divv\mathbf{u}_t\mathbb{I}\right) \nonumber \\
& = -\nabla\partial_{t}\rho^{\gamma}
+\divv\left(2\nu_{t}\mathfrak{D}(\mathbf{u})
+\lambda_{t}\divv\mathbf{u}\mathbb{I}\right)
-\rho_{t}\left(\mathbf{u}_{t}+\mathbf{u}\cdot\nabla\mathbf{u}\right)
\nonumber \\
& \quad -\rho\mathbf{u}_{t}\cdot\nabla\mathbf{u}
-\mathbf{H}_{t}\cdot\nabla\mathbf{H}-\mathbf{H}\cdot\nabla\mathbf{H}_{t}
-\frac{1}{2}\nabla\partial_{t}|\mathbf{H}|^2,
\end{align}
where $\mathbb{I}$ is the identity matrix. Multiplying \eqref{3.21} by $\mathbf{u}_t$ and integrating over $\Omega$, we have
\begin{align}\label{3.22}
& \frac{1}{2}\frac{d}{dt}\int\rho|\mathbf{u}_{t}|^2dx
+\int\left(\nu|\nabla\mathbf{u}_{t}|^2+\mu|\divv\mathbf{u}_{t}|^2
\right)dx\nonumber\\
& = \int\gamma\rho^{\gamma-1}\rho_{t}\divv\mathbf{u}_{t}dx
-\int\left(2\nu'\rho_{t}\mathfrak{\mathbf{u}}
+\lambda'\rho_{t}\divv\mathbf{u}_{t}\mathbb{I}\right)\divv\mathbf{u}_{t}dx
\nonumber\\
& \quad
-\int\rho_{t}|\mathbf{u}_{t}|^2dx
-\int\rho_{t}\mathbf{u}\cdot\nabla\mathbf{u}\cdot\mathbf{u}_{t}dx
-\int\rho\mathbf{u}_{t}\cdot\nabla\mathbf{u}\cdot\mathbf{u}_{t}dx
\nonumber\\
& \quad
-\int\mathbf{H}_{t}\cdot\nabla\mathbf{H}\cdot\mathbf{u}_{t}dx
-\int\mathbf{H}\cdot\nabla\mathbf{H}_{t}\cdot\mathbf{u}_{t}dx
-\frac{1}{2}\int\nabla\partial_{t}|\mathbf{H}|^2\cdot\mathbf{u}_{t}dx
\triangleq \sum_{k=1}^{8}J_k.
\end{align}
We can bound each term $J_k$ on the right-hand side of \eqref{3.22} as follows.
\begin{align*}
|J_{1}| \leq & C\|\rho^{\gamma-1}\|_{L^\infty}\|\rho_t\|_{L^2}
\|\nabla\mathbf{u}_{t}\|_{L^2}
\leq \varepsilon\|\nabla\mathbf{u}_{t}\|_{L^2}^2
+C(\varepsilon)\|\rho_t\|_{L^2}^2; \\
|J_{2}| \leq & C\|\rho_t\|_{L^2}\|\nabla\mathbf{u}\|_{L^\infty}
\|\nabla\mathbf{u}_{t}\|_{L^2}
\leq \varepsilon\|\nabla\mathbf{u}_{t}\|_{L^2}^2
+C(\varepsilon)\|\nabla\mathbf{u}\|_{L^\infty}^2\|\rho_t\|_{L^2}^2;\\
|J_{3}|= & \left|\int\divv(\rho\mathbf{u})|\mathbf{u}_t|^2dx\right|
=\left|-\int\rho\mathbf{u}\cdot\nabla|\mathbf{u}_t|^2dx\right|
\quad (\text{by}\ \eqref{1.1}_1) \\
\leq & C\|\mathbf{u}\|_{L^\infty}\|\sqrt{\rho}\|_{L^\infty}
\|\sqrt{\rho}\mathbf{u}_{t}\|_{L^2}\|\nabla\mathbf{u}_{t}\|_{L^2} \\ \leq & \varepsilon\|\nabla\mathbf{u}_{t}\|_{L^2}^2
+C(\varepsilon)\|\nabla\mathbf{u}\|_{L^\infty}^2
\|\sqrt{\rho}\mathbf{u}_t\|_{L^2}^2;\\
|J_{4}|\leq &\|\mathbf{u}\|_{L^\infty}
\|\nabla\mathbf{u}\|_{L^\infty}\|\rho^{-\frac{1}{2}}\|_{L^\infty}\|\rho_t\|_{L^2}
\|\sqrt{\rho}\mathbf{u}_{t}\|_{L^2}
\leq  C\|\nabla\mathbf{u}\|_{L^\infty}^2\|\rho_t\|_{L^2}
\|\sqrt{\rho}\mathbf{u}_{t}\|_{L^2} \\
\leq & C\|\nabla\mathbf{u}\|_{L^\infty}^2
\left(\|\rho_t\|_{L^2}^2+\|\sqrt{\rho}\mathbf{u}_{t}\|_{L^2}^2\right);\\
|J_{5}|\leq & \|\nabla\mathbf{u}\|_{L^\infty}
\|\sqrt{\rho}\mathbf{u}_{t}\|_{L^2}^2
\leq C\left(1+\|\nabla\mathbf{u}\|_{L^\infty}^2\right)
\|\sqrt{\rho}\mathbf{u}_{t}\|_{L^2}^2;
\\
\sum_{k=6}^{8}|J_{k}|\leq & C\int|\mathbf{H}||\mathbf{H}_t||\nabla\mathbf{u}_t|dx
\leq \varepsilon\|\nabla\mathbf{u}_{t}\|_{L^2}^2
+C(\varepsilon)\||\mathbf{H}||\mathbf{H}_t|\|_{L^2}^2 \\
\leq & \varepsilon\|\nabla\mathbf{u}_{t}\|_{L^2}^2
+C(\varepsilon)\|\mathbf{H}\|_{L^6}^2
\|\mathbf{H}_t\|_{L^2}\|\mathbf{H}_t\|_{L^6}
\\ \leq & \varepsilon\|\nabla\mathbf{u}_{t}\|_{L^2}^2
+C(\varepsilon,\delta)
\|\mathbf{H}_t\|_{L^2}^2+\delta\|\nabla\mathbf{H}_t\|_{L^6}^2
\end{align*}
due to Lemma \ref{lem22}, \eqref{3.2}, \eqref{a}, and \eqref{3.5}. Hence, choosing $\varepsilon$ suitably small and inserting the above estimates into \eqref{3.22}, we derive from \eqref{1.2} that
\begin{align}\label{3.23}
& \frac{d}{dt}\int\rho|\mathbf{u}_{t}|^2dx
+C_5\int|\nabla\mathbf{u}_{t}|^2dx\nonumber\\
& \leq C\left(1+\|\nabla\mathbf{u}\|_{L^\infty}^2\right)
\left(\|\rho_{t}\|_{L^2}^2+
\|\sqrt{\rho}\mathbf{u}_{t}\|_{L^2}^2+\|\mathbf{H}_{t}\|_{L^2}^2\right)
+2\delta\|\nabla\mathbf{H}_t\|_{L^2}^2
\end{align}
for some $C_5>0$.

Next, differentiating \eqref{1.1}$_3$ with respect to $t$ leads to
\begin{equation}\label{4.1}
\mathbf{H}_{tt}-\mathbf{H}_{t}\cdot\nabla\mathbf{u}
-\mathbf{H}\cdot\nabla\mathbf{u}_{t}
+\mathbf{u}_{t}\cdot\nabla\mathbf{H}
+\mathbf{u}\cdot\nabla\mathbf{H}_{t}
+\mathbf{H}_{t}\divv\mathbf{u}
+\mathbf{H}\divv\mathbf{u}_{t}=\beta\Delta\mathbf{H}_{t}.
\end{equation}
Multiplying \eqref{4.1} by $\mathbf{H}_t$ and integrating the resulting equations by parts show that
\begin{align}\label{4.2}
& \frac{1}{2}\frac{d}{dt}\int|\mathbf{H}_{t}|^2dx
+\beta\int|\nabla\mathbf{H}_{t}|^2dx \nonumber \\
& = -\int\mathbf{H}\cdot\nabla\mathbf{H}_{t}\cdot\mathbf{u}_{t}dx
+\int\mathbf{u}_{t}\cdot\nabla\mathbf{H}_{t}\cdot\mathbf{H}dx
-\int\mathbf{H}_{t}\cdot\nabla\mathbf{H}_{t}\cdot\mathbf{u}dx
+\int\mathbf{u}\cdot\nabla\mathbf{H}_{t}\cdot\mathbf{H}_{t}dx
\nonumber \\
& \leq \frac{\beta}{4}\|\nabla\mathbf{H}_{t}\|_{L^2}^2
+C\||\mathbf{u}_{t}||\mathbf{H}|\|_{L^2}^2
+C\||\mathbf{u}||\mathbf{H}_{t}|\|_{L^2}^2
\nonumber \\
& \leq \frac{\beta}{4}\|\nabla\mathbf{H}_{t}\|_{L^2}^2
+C\|\mathbf{u}\|_{L^6}^2\|\mathbf{H}_{t}\|_{L^3}^2
+C\|\mathbf{H}\|_{L^6}^2\|\mathbf{u}_{t}\|_{L^3}^2
\nonumber \\
& \leq \frac{\beta}{4}\|\nabla\mathbf{H}_{t}\|_{L^2}^2
+C\|\mathbf{u}\|_{L^6}^2\|\mathbf{H}_{t}\|_{L^2}\|\mathbf{H}_{t}\|_{L^6}
+C\|\mathbf{H}\|_{L^6}^2\|\mathbf{u}_{t}\|_{L^2}\|\mathbf{u}_{t}\|_{L^6}
\nonumber \\
& \leq \frac{\beta}{2}\|\nabla\mathbf{H}_{t}\|_{L^2}^2
+C\|\mathbf{H}_{t}\|_{L^2}^2
+C\|\nabla\mathbf{u}_{t}\|_{L^2}^2
\end{align}
owing to \eqref{3.5}, \eqref{3.7}, and Lemma \ref{lem22}. Thus we get
\begin{align}\label{4.3}
\frac{d}{dt}\int|\mathbf{H}_{t}|^2dx
+\beta\int|\nabla\mathbf{H}_{t}|^2dx
\leq C\|\mathbf{H}_{t}\|_{L^2}^2
+C\|\nabla\mathbf{u}_{t}\|_{L^2}^2.
\end{align}
Adding \eqref{3.23} multiplied by $C_6$ suitably large, and choosing $\delta$ small enough such that $2C_6\delta<\beta$, to \eqref{4.3}, we deduce from \eqref{x} and \eqref{3.7} that
\begin{align}\label{4.4}
& \frac{d}{dt}\int\left(\rho|\mathbf{u}_{t}|^2+|\mathbf{H}_{t}|^2\right)dx
+C_7\int\left(|\nabla\mathbf{u}_{t}|^2+|\nabla\mathbf{H}_{t}|^2\right)dx\nonumber\\
& \leq C \left(1+\|\nabla\mathbf{u}\|_{L^\infty}^2\right)
\left(\|\rho_{t}\|_{L^2}^2+
\|\sqrt{\rho}\mathbf{u}_{t}\|_{L^2}^2+\|\mathbf{H}_{t}\|_{L^2}^2\right)
\nonumber\\
& \leq C \left(1+\|\nabla\mathbf{u}\|_{L^\infty}^2\right)
\left(1+
\|\sqrt{\rho}\mathbf{u}_{t}\|_{L^2}^2+\|\mathbf{H}_{t}\|_{L^2}^2\right)
\end{align}
for some $C_7>0$.
Then we obtain from the Gronwall inequality and \eqref{3.1} that
\begin{align}\label{4.6}
\sup_{0\leq t\leq T}\int\left(\rho|\mathbf{u}_{t}|^2+|\mathbf{H}_{t}|^2\right)dx
+\int_{0}^{T}\int\left(|\nabla\mathbf{u}_{t}|^2
+|\nabla\mathbf{H}_{t}|^2\right)dxdt \leq C.
\end{align}


Now we are going to bound $\|\rho\|_{H^2}$. To this end, we rewrite  \eqref{x2} in the following
\begin{align}\label{z3}
\begin{cases}
\mathfrak{L}\mathbf{u}=\frac{\mathbf{F}}{\nu},\ \ x\in\Omega,\\
\mathbf{u}=\mathbf{0},\ \ x\in\partial\Omega.
\end{cases}
\end{align}
Here
\begin{equation*}
\mathfrak{L}\mathbf{u}\triangleq\Delta\mathbf{u}
+\left(k+1\right)\nabla\divv\mathbf{u}
\end{equation*}
is the Lam{\'e} operator and $\mathbf{F}$ is defined by \eqref{x3}. By standard elliptic estimates, we infer from \eqref{z3} that
\begin{align}\label{z4}
\|\mathbf{u}\|_{H^3}^2
& \leq C\|\rho\mathbf{u}_{t}\|_{H^1}^2
+C\|\rho|\mathbf{u}||\nabla\mathbf{u}|\|_{H^1}^2
+C\|\nabla\rho^{\gamma}\|_{H^1}^2
+C\||\nabla\nu||\nabla\mathbf{u}|\|_{H^1}^2 \nonumber\\
& \quad
+C\||\nabla\lambda||\nabla\mathbf{u}|\|_{H^1}^2
+C\||\mathbf{H}||\nabla\mathbf{H}|\|_{H^1}^2
+C\|\mathbf{u}\|_{H^1}^2 \triangleq\sum_{j=1}^{7}\bar{I_j}.
\end{align}
By Lemma \ref{lem22}, \eqref{3.2}, \eqref{3.7}, and \eqref{4.6}, one has
\begin{equation*}
\begin{split}
\bar{I_1} & \leq C\|\rho\|_{L^\infty}\|\sqrt{\rho}\mathbf{u}_{t}\|_{L^2}^2
+C\||\nabla\rho||\mathbf{u}_t|\|_{L^2}^2
+C\|\rho\|_{L^\infty}^2\|\nabla\mathbf{u}_{t}\|_{L^2}^2 \\
& \leq C+C\|\nabla\rho\|_{H^1}^2\|\mathbf{u}_{t}\|_{L^2}
\|\mathbf{u}_{t}\|_{L^6}
+C\|\nabla\mathbf{u}_{t}\|_{L^2}^2 \\
& \leq C\left(1+\|\nabla\mathbf{u}_t\|_{L^2}^2\right)
\left(1+\|\nabla^2\rho\|_{L^2}^2\right).
\end{split}
\end{equation*}
Similarly, one gets
\begin{equation*}
\begin{split}
\bar{I_2} & \leq C\|\rho\|_{L^\infty}\|\mathbf{u}\|_{L^\infty}\|\nabla\mathbf{u}\|_{L^2}^2
+C\||\nabla\rho||\mathbf{u}||\nabla\mathbf{u}|\|_{L^2}^2
+C\|\rho|\nabla\mathbf{u}|^2\|_{L^2}^2
+C\|\rho|\mathbf{u}||\nabla^2\mathbf{u}|\|_{L^2}^2 \\
& \leq C\|\mathbf{u}\|_{L^\infty}
+C\|\nabla\rho\|_{L^6}^2\|\mathbf{u}\|_{L^6}^2\|\nabla\mathbf{u}\|_{L^6}^2
+C\|\rho\|_{L^\infty}^2\|\nabla\mathbf{u}\|_{L^\infty}^2
\|\nabla\mathbf{u}\|_{L^2}^2
+C\|\rho\|_{L^\infty}^2\|\mathbf{u}\|_{L^6}^2
\|\nabla^2\mathbf{u}\|_{L^3}^2
\\
& \leq C\|\nabla\mathbf{u}\|_{L^\infty}
+C\left(1+\|\nabla^2\rho\|_{L^2}^2\right)
\left(1+\|\nabla^2\mathbf{u}\|_{L^2}^2\right)
+C\|\nabla\mathbf{u}\|_{L^\infty}^2
+C\|\nabla^2\mathbf{u}\|_{L^2}\|\nabla^2\mathbf{u}\|_{L^6}
\\
& \leq C\left(1+\|\nabla\mathbf{u}\|_{L^\infty}^2+\|\nabla^2\mathbf{u}\|_{L^2}^2\right)
\left(1+\|\nabla^2\rho\|_{L^2}^2\right)
+\frac{1}{8}\|\mathbf{u}\|_{H^3}^2.
\end{split}
\end{equation*}
By \eqref{3.2}, \eqref{3.3}, and \eqref{3.7}, we have
\begin{equation*}
\bar{I_3}+\bar{I_7} \leq C+C\|\nabla^2\rho\|_{L^2}^2.
\end{equation*}
We temporarily claim that for any $2< p\leq6$, we have obtained the following estimates on $\|\nabla\rho\|_{L^p}$,
\begin{equation}\label{9.1}
\sup_{0\leq t\leq T}\|\nabla\rho\|_{L^p}\leq C,
\end{equation}
whose proof can be found in Appendix of this paper.
Then by \eqref{z}, \eqref{3.7}, and \eqref{9.1}, we arrive at
\begin{equation*}
\begin{split}
\bar{I_4} & \leq C\||\nu'||\nabla\rho||\nabla\mathbf{u}|\|_{L^2}^2
+C\||\nu''||\nabla\rho|^2|\nabla\mathbf{u}|\|_{L^2}^2
+C\||\nu'||\nabla^2\rho||\nabla\mathbf{u}|\|_{L^2}^2
+C\||\nu'||\nabla\rho||\nabla^2\mathbf{u}|\|_{L^2}^2
\\ & \leq C\|\nabla\mathbf{u}\|_{L^\infty}^2\|\nabla\rho\|_{L^2}^2
+C\|\nabla\mathbf{u}\|_{L^\infty}^2\|\nabla\rho\|_{L^4}^2
+C\|\nabla\mathbf{u}\|_{L^\infty}^2
\|\nabla^2\rho\|_{L^2}^2
+C\|\nabla\rho\|_{L^6}^2
\|\nabla^2\mathbf{u}\|_{L^2}\|\nabla^2\mathbf{u}\|_{L^6}
\\
& \leq C\left(1+\|\nabla\mathbf{u}\|_{L^\infty}^2
+\|\nabla^2\mathbf{u}\|_{L^2}^2\right)
\left(1+\|\nabla^2\rho\|_{L^2}^2\right)
+\frac{1}{8}\|\mathbf{u}\|_{H^3}^2.
\end{split}
\end{equation*}
Similarly to $\bar{I_4}$, we derive
\begin{equation*}
\bar{I_5} \leq C\left(1+\|\nabla\mathbf{u}\|_{L^\infty}^2
+\|\nabla^2\mathbf{u}\|_{L^2}^2\right)
\left(1+\|\nabla^2\rho\|_{L^2}^2\right)
+\frac{1}{8}\|\mathbf{u}\|_{H^3}^2.
\end{equation*}
From \eqref{3.5} and \eqref{3.7}, we get
\begin{equation*}
\begin{split}
\bar{I_6} & \leq C\||\mathbf{H}||\nabla\mathbf{H}|\|_{L^2}^2
+C\||\nabla\mathbf{H}|^2\|_{L^2}^2
+C\||\mathbf{H}||\nabla^2\mathbf{H}|\|_{L^2}^2 \\
& \leq C\|\mathbf{H}\|_{L^6}^2\|\nabla\mathbf{H}\|_{L^2}
\|\nabla\mathbf{H}\|_{L^6}
+C\|\nabla\mathbf{H}\|_{L^4}^2
+C\|\mathbf{H}\|_{L^6}^2\|\nabla^2\mathbf{H}\|_{L^2}
\|\nabla^2\mathbf{H}\|_{L^6}
\\
& \leq C+C\|\nabla^2\mathbf{H}\|_{L^2}^2+\delta\|\mathbf{H}\|_{H^3}^2.
\end{split}
\end{equation*}
In order to control $\|\mathbf{H}\|_{H^3}^2$, by virtue of standard regularity theory of elliptic system, we obtain from \eqref{1.1}$_3$, \eqref{4.6}, \eqref{3.3}, \eqref{3.5}, \eqref{3.7}, and Lemma \ref{lem22} that
\begin{equation*}\label{8.2}
\begin{split}
\|\mathbf{H}\|_{H^3}^2
& \leq C\|\mathbf{H}_t\|_{H^1}^2
+C\||\mathbf{H}||\nabla\mathbf{u}|\|_{H^1}^2
+C\||\mathbf{u}||\nabla\mathbf{H}|\|_{H^1}^2+C\|\mathbf{H}\|_{H^1}^2 \nonumber\\
& \leq C+C\|\nabla\mathbf{H}_t\|_{L^2}^2
+C\|\nabla\mathbf{u}\|_{L^\infty}^2\|\mathbf{H}\|_{L^2}^2
+C\|\nabla\mathbf{u}\|_{L^\infty}^2\|\nabla\mathbf{H}\|_{L^2}^2
\nonumber\\
& \quad +C\|\mathbf{H}\|_{L^6}^2\|\nabla^2\mathbf{u}\|_{L^2}
\|\nabla^2\mathbf{u}\|_{L^6}
+C\|\mathbf{u}\|_{L^\infty}^2\|\nabla\mathbf{H}\|_{L^2}^2
+C\|\mathbf{u}\|_{L^6}^2\|\nabla^2\mathbf{H}\|_{L^2}\|\nabla^2\mathbf{H}\|_{L^6}
\nonumber\\
& \leq C\left(1+\|\nabla\mathbf{H}_t\|_{L^2}^2
+\|\nabla\mathbf{u}\|_{L^\infty}^2
+\|\nabla^2\mathbf{u}\|_{L^2}^2
+\|\nabla^2\mathbf{H}\|_{L^2}^2\right)
+\frac{1}{2}\|\mathbf{H}\|_{H^3}^2+C\|\mathbf{u}\|_{H^3}^2.
\end{split}
\end{equation*}
Therefore, one has
\begin{align}\label{8.2}
\|\mathbf{H}\|_{H^3}^2\leq
C\left(1+\|\nabla\mathbf{H}_t\|_{L^2}^2
+\|\nabla\mathbf{u}\|_{L^\infty}^2
+\|\nabla^2\mathbf{u}\|_{L^2}^2
+\|\nabla^2\mathbf{H}\|_{L^2}^2\right)
+C\|\mathbf{u}\|_{H^3}^2.
\end{align}
Inserting the above estimates on $\bar{I_j}$ into \eqref{z4} and choosing $\delta$ suitably small and together with \eqref{8.2}, we find
\begin{align}\label{z10}
\|\mathbf{u}\|_{H^3}^2
\leq CC_2(t)\left(1+\|\nabla^2\rho\|_{L^2}^2\right),
\end{align}
where
\begin{equation}\label{t}
C_2(t)\triangleq 1+\|\nabla\mathbf{u}\|_{L^\infty}^2
+\|\nabla\mathbf{u}_t\|_{L^2}^2
+\|\nabla\mathbf{H}_t\|_{L^2}^2
+\|\nabla^2\mathbf{u}\|_{L^2}^2
+\|\nabla^2\mathbf{H}\|_{L^2}^2\in L^{1}(0,T)
\end{equation}
due to \eqref{3.1}, \eqref{4.6}, \eqref{8.1}, and \eqref{3.7}. Moreover, taking the operator $\frac{\partial^2}{\partial x_j\partial x_k}$ on both side of \eqref{1.1}$_1$, and then multiplying the resulting equations by $\frac{\partial^2\rho}{\partial x_j\partial x_k}$, we deduce from Lemma \ref{lem22}, \eqref{3.2}, \eqref{3.7}, \eqref{9.1}, and \eqref{z10} that
\begin{align}\label{10}
\frac{d}{dt}\int|\nabla^2\rho|^2dx
& \leq C\int|\nabla^2\rho|\left(|\nabla^2\rho||\nabla\mathbf{u}|
+|\nabla\rho||\nabla^2\mathbf{u}|
+|\rho||\nabla^3\mathbf{u}|\right)dx
 \nonumber \\
& \leq C\|\nabla\mathbf{u}\|_{L^\infty}\|\nabla^2\rho\|_{L^2}^2
+C\|\nabla^2\rho\|_{L^2}\||\nabla\rho||\nabla^2\mathbf{u}|\|_{L^2}
+C\|\rho\|_{L^\infty}\|\nabla^3\mathbf{u}\|_{L^2}\|\nabla^2\rho\|_{L^2}
 \nonumber \\
& \leq C\left(1+\|\nabla\mathbf{u}\|_{L^\infty}^2\right)
\|\nabla^2\rho\|_{L^2}^2+C\|\nabla\rho\|_{L^6}^{2}
\|\nabla^2\mathbf{u}\|_{L^2}\|\nabla^2\mathbf{u}\|_{L^6}
+C\|\mathbf{u}\|_{H^3}^2
 \nonumber \\
& \leq C\left(1+\|\nabla\mathbf{u}\|_{L^\infty}^2
+\|\nabla^2\mathbf{u}\|_{L^2}^2\right)
\left(1+\|\nabla^2\rho\|_{L^2}^2\right)
+C\|\mathbf{u}\|_{H^3}^2
 \nonumber \\
& \leq CC_2(t)\left(1+\|\nabla^2\rho\|_{L^2}^2\right)
\end{align}
with $C_{2}(t)$ as in \eqref{t}, where we have used
\begin{equation*}
-\int\partial_{ijk}\rho u^{i}\partial_{jk}\rho dx
=\int\partial_{jk}\rho\partial_{i}u^{i}\partial_{jk}\rho dx
+\int\partial_{jk}\rho u^{i}\partial_{ijk}\rho dx,
\end{equation*}
that is,
\begin{equation*}
-\int\partial_{ijk}\rho u^{i}\partial_{jk}\rho dx
=\frac{1}{2}\int\partial_{jk}\rho\partial_{i}u^{i}\partial_{jk}\rho dx
\end{equation*}
in the derivation of the first inequality.
Hence we conclude by Gronwall's inequality and \eqref{t} that
\begin{equation*}
\sup_{0\leq t\leq T}\|\nabla^2\rho\|_{L^2}\leq C,
\end{equation*}
which along with \eqref{3.2} and \eqref{3.7} leads to
\begin{equation}\label{z2}
\sup_{0\leq t\leq T}\|\rho\|_{H^2}\leq C.
\end{equation}

We turn to estimate $\|\mathbf{u}\|_{H^2}$.
It follows from \eqref{x2}, \eqref{x3}, and Lemma \ref{lem23} that
\begin{align}\label{x4}
\|\nabla^2\mathbf{u}\|_{L^2}
& \leq C\|\rho\partial_{t}\mathbf{u}
+\rho\mathbf{u}\cdot\nabla\mathbf{u}
+\nabla\rho^\gamma
-\nabla\nu\cdot\nabla\mathbf{u}
-\nabla\mathbf{u}\cdot\nabla\nu
-\divv\mathbf{u}\nabla\lambda
-(\nabla\times\mathbf{H})\times\mathbf{H}\|_{L^2}
\nonumber\\
& \leq C\|\sqrt{\rho}\|_{L^\infty}\|\sqrt{\rho}\partial_{t}\mathbf{u}\|_{L^2}
+C\|\rho\|_{L^\infty}\||\mathbf{u}||\nabla\mathbf{u}|\|_{L^2}
+C\|\rho^{\gamma-1}\|_{L^\infty}\|\nabla\rho\|_{L^2} \nonumber\\
& \quad
+C\|\nu'\|_{L^\infty}\||\nabla\rho||\nabla\mathbf{u}|\|_{L^2}
+C\|\lambda'\|_{L^\infty}\||\nabla\rho||\nabla\mathbf{u}|\|_{L^2}
+C\|\nabla\mathbf{H}\|_{L^2}^{\frac{1}{2}}\|\nabla\mathbf{H}\|_{H^1}^{\frac{1}{2}}
 \nonumber\\
& \triangleq\sum_{j=1}^{6}K_j.
\end{align}
By \eqref{3.2}, \eqref{4.6}, and \eqref{3.7}, one has
\begin{equation*}
K_1+K_3\leq C.
\end{equation*}
By Lemma \ref{lem22}, we get from \eqref{3.7} that
\begin{equation*}
\begin{split}
K_2 \leq C\|\mathbf{u}\|_{L^6}\|\nabla\mathbf{u}\|_{L^2}^{\frac{1}{2}}
\|\nabla\mathbf{u}\|_{L^6}^{\frac{1}{2}}
\leq C\|\nabla\mathbf{u}\|_{L^2}^2
+C\|\nabla\mathbf{u}\|_{L^2}^{\frac{3}{2}}
\|\nabla^2\mathbf{u}\|_{L^2}^{\frac{1}{2}}
\leq C+\frac{1}{4}\|\nabla^2\mathbf{u}\|_{L^2}.
\end{split}
\end{equation*}
We obtain from \eqref{z}, Lemma \ref{lem22}, \eqref{3.7}, and \eqref{z2} that
\begin{equation*}
\begin{split}
K_4+K_5 & \leq C\|\nabla\rho\|_{L^6}\|\nabla\mathbf{u}\|_{L^2}^{\frac{1}{2}}
\|\nabla\mathbf{u}\|_{L^6}^{\frac{1}{2}} \\
& \leq C\|\rho\|_{H^2}\|\nabla\mathbf{u}\|_{L^2}+
C\|\rho\|_{H^2}\|\nabla\mathbf{u}\|_{L^2}^{\frac{1}{2}}
\|\nabla^2\mathbf{u}\|_{L^2}^{\frac{1}{2}}
\\ & \leq C+\frac{1}{4}\|\nabla^2\mathbf{u}\|_{L^2}.
\end{split}
\end{equation*}
Sobolev's inequality and \eqref{3.7} give that
\begin{equation*}
K_6  \leq C\|\nabla\mathbf{H}\|_{L^2}+C\|\nabla\mathbf{H}\|_{L^2}^{\frac{1}{2}}
\|\nabla^2\mathbf{H}\|_{L^2}^{\frac{1}{2}} \leq C+\varepsilon\|\mathbf{H}\|_{H^2}.
\end{equation*}
By standard regularity theory of elliptic system, it follows from \eqref{1.1}$_3$, \eqref{3.3}, \eqref{3.5}, \eqref{3.7}, \eqref{4.6}, and Lemma \ref{lem22} that
\begin{align}\label{4.7}
\|\mathbf{H}\|_{H^2}^2
& \leq C\|\mathbf{H}_t\|_{L^2}^2
+C\||\mathbf{H}||\nabla\mathbf{u}|\|_{L^2}^2
+C\||\mathbf{u}||\nabla\mathbf{H}|\|_{L^2}^2+C\|\mathbf{H}\|_{L^2}^2 \nonumber\\
& \leq C+C\|\nabla\mathbf{u}\|_{L^4}^2+C\|\nabla\mathbf{H}\|_{L^4}^2
\nonumber\\
& \leq C+C\|\nabla^2\mathbf{u}\|_{L^2}^2
+\frac{1}{2}\|\nabla^2\mathbf{H}\|_{L^2}^2,
\end{align}
which leads to
\begin{equation}\label{4.8}
  \|\mathbf{H}\|_{H^2}^2\leq C+C\|\nabla^2\mathbf{u}\|_{L^2}^2,
\end{equation}
which along with the above estimates on $K_j$, we conclude by choosing $\varepsilon$ suitably small that
\begin{equation*}
\sup_{0\leq t\leq T}\|\nabla^2\mathbf{u}\|_{L^2}\leq C,
\end{equation*}
which combined with \eqref{3.3}, \eqref{3.2}, and \eqref{3.7} implies that
\begin{equation}\label{4.12}
\sup_{0\leq t\leq T}\|\mathbf{u}\|_{H^2}\leq C.
\end{equation}

Finally, by Sobolev's inequality, \eqref{4.8}, and \eqref{4.12}, we have
\begin{equation}\label{4.13}
\sup_{0\leq t\leq T}\|\mathbf{H}\|_{H^2}
\leq C+C\sup_{0\leq t\leq T}\|\mathbf{u}\|_{H^2}\leq C.
\end{equation}

Then the desired \eqref{123} follows from \eqref{z2}, \eqref{4.12}, and \eqref{4.13}. This completes the proof of Lemma \ref{lem35}.
\hfill $\Box$

\begin{lemma}\label{lem36}
Under the condition \eqref{3.1}, it holds that for any $T\in[0,T^*)$,
\begin{equation}\label{l}
\sup_{0\leq t\leq T}\left(\|\mathbf{u}\|_{H^3}
+\|\rho\|_{H^3}+\|\mathbf{H}\|_{H^3}\right) \leq C.
\end{equation}
\end{lemma}
{\it Proof.}
The proof of $L^{\infty}(0,T;H^3)$-norm of strong solution $(\rho,\mathbf{u},\mathbf{H})$ follows exactly the same as the derivation of
$L^{\infty}(0,T;H^2)$-norm of $(\rho,\mathbf{u},\mathbf{H})$, hence we sketch it for simplicity.

First, we shall show that
\begin{align}\label{zx}
\sup_{0\leq t\leq T}\int\left(\rho|\nabla\mathbf{u}_{t}|^2+|\nabla\mathbf{H}_{t}|^2\right)dx
+\int_{0}^{T}\int\left(|\nabla^2\mathbf{u}_{t}|^2
+|\nabla^2\mathbf{H}_{t}|^2\right)dxdt \leq C.
\end{align}
On one hand, taking the operator $\frac{\partial}{\partial t}\frac{\partial}{\partial x_i}$ on both side of \eqref{1.1}$_3$, and then multiplying the resulting equations by $\partial_{i}\mathbf{H}_t$, we obtain that
\begin{equation*}\label{zx2}
\begin{split}
& \frac{1}{2}\frac{d}{dt}\int|\nabla\mathbf{H}_{t}|^2dx
+\beta\int|\nabla^2\mathbf{H}_{t}|^2dx \nonumber \\
& \leq C\int|\nabla\mathbf{H}_{t}|^2|\nabla\mathbf{u}|dx
+C\int|\mathbf{H}_t||\nabla\mathbf{H}_{t}||\nabla^2\mathbf{u}|dx
+C\int|\nabla\mathbf{H}||\nabla\mathbf{H}_{t}||\nabla\mathbf{u}_{t}|dx
\nonumber \\
& \quad
+C\int|\mathbf{H}||\nabla\mathbf{H}_{t}||\nabla^2\mathbf{u}_{t}|dx
+C\int|\mathbf{u}_t||\nabla\mathbf{H}_{t}||\nabla^2\mathbf{H}|dx
+C\int|\mathbf{u}||\nabla\mathbf{H}_{t}||\nabla^2\mathbf{H}_{t}|dx
\nonumber \\
& \leq C\|\nabla\mathbf{u}\|_{L^\infty}\|\nabla\mathbf{H}_{t}\|_{L^2}^2
+C\|\mathbf{H}_t\|_{L^\infty}\|\nabla\mathbf{H}_{t}\|_{L^2}
\|\nabla^2\mathbf{u}\|_{L^2}
+C\|\nabla\mathbf{H}\|_{L^6}\|\nabla\mathbf{H}_{t}\|_{L^3}
\|\nabla\mathbf{u}_t\|_{L^2}
\nonumber \\
& \quad
+C\|\mathbf{H}\|_{L^\infty}\|\nabla\mathbf{H}_{t}\|_{L^2}
\|\nabla^2\mathbf{u}_{t}\|_{L^2}
+C\|\mathbf{u}_t\|_{L^\infty}\|\nabla\mathbf{H}_{t}\|_{L^2}
\|\nabla^2\mathbf{H}\|_{L^2}
+C\|\mathbf{u}\|_{L^\infty}\|\nabla\mathbf{H}_{t}\|_{L^2}
\|\nabla^2\mathbf{H}_t\|_{L^2}
\nonumber \\
& \leq \frac{\beta}{2}\|\nabla^2\mathbf{H}_{t}\|_{L^2}^2
+C\left(1+\|\nabla\mathbf{u}\|_{L^\infty}^2+\|\nabla\mathbf{u}_{t}\|_{L^2}^2\right)
\left(1+\|\nabla\mathbf{H}_{t}\|_{L^2}^2\right)
+\varepsilon\|\nabla^2\mathbf{u}_{t}\|_{L^2}^2,
\end{split}
\end{equation*}
where one has used H{\"o}lder's inequality, \eqref{123}, and \eqref{4.6}. Then we have
\begin{align}\label{zx2}
\frac{d}{dt}\int|\nabla\mathbf{H}_{t}|^2dx
+\beta\int|\nabla^2\mathbf{H}_{t}|^2dx
\leq CC_{3}(t)\left(1+\|\nabla\mathbf{H}_{t}\|_{L^2}^2\right)
+\varepsilon\|\nabla^2\mathbf{u}_{t}\|_{L^2}^2,
\end{align}
where
\begin{equation*}
C_{3}(t)\triangleq
1+\|\nabla\mathbf{u}\|_{L^\infty}^2
+\|\nabla\mathbf{u}_{t}\|_{L^2}^2\in L^{1}(0,T)
\end{equation*}
due to \eqref{3.1} and \eqref{4.6}.
On the other hand, taking the operator $\frac{\partial}{\partial t}\frac{\partial}{\partial x_i}$ on both side of \eqref{1.1}$_2$, and then multiplying the resulting equations by $\partial_{i}\mathbf{u}_t$, we can deduce that
\begin{align}\label{5.1}
& \frac{d}{dt}\int\rho|\nabla\mathbf{u}_{t}|^2dx
+C\int|\nabla^2\mathbf{u}_{t}|^2dx
\nonumber \\
& \leq C\left(1+\|\nabla\mathbf{u}\|_{L^\infty}^2\right)
\left(1+\|\nabla\rho_{t}\|_{L^2}^2
+\|\sqrt{\rho}\nabla\mathbf{u}_{t}\|_{L^2}^2
+\|\nabla\mathbf{H}_{t}\|_{L^2}^2\right)
\nonumber \\
& \quad
+C\|\nabla^2\mathbf{H}_{t}\|_{L^2}^2+\delta\|\nabla^2\mathbf{u}_{t}\|_{L^2}^2,
\end{align}
where one has used H{\"o}lder's inequality and \eqref{123}. Hence, choosing $\delta$ in \eqref{5.1} suitably small, we get
\begin{align}\label{5.2}
& \frac{d}{dt}\int\rho|\nabla\mathbf{u}_{t}|^2dx
+\widetilde{C_1}\int|\nabla^2\mathbf{u}_{t}|^2dx
\nonumber \\
& \leq C\left(1+\|\nabla\mathbf{u}\|_{L^\infty}^2\right)
\left(1+\|\nabla\rho_{t}\|_{L^2}^2
+\|\sqrt{\rho}\nabla\mathbf{u}_{t}\|_{L^2}^2
+\|\nabla\mathbf{H}_{t}\|_{L^2}^2\right)
+C\|\nabla^2\mathbf{H}_{t}\|_{L^2}^2
\end{align}
for some $\widetilde{C_1}>0$.
Then, adding \eqref{zx2} multiplied by $\widetilde{C_2}$ suitably large, and choosing $\varepsilon$ in \eqref{zx2} small enough to \eqref{5.2}, we derive that
\begin{align}\label{5.3}
& \frac{d}{dt}\int\left(\rho|\nabla\mathbf{u}_{t}|^2
+|\nabla\mathbf{H}_{t}|^2\right)dx
+C_8\int\left(|\nabla^2\mathbf{u}_{t}|^2+
|\nabla^2\mathbf{H}_{t}|^2\right)dx\nonumber\\
& \leq CC_{3}(t)\left(1+\|\nabla\rho_{t}\|_{L^2}^2
+\|\sqrt{\rho}\nabla\mathbf{u}_{t}\|_{L^2}^2
+\|\nabla\mathbf{H}_{t}\|_{L^2}^2\right)
\end{align}
for some $C_8>0$. From the continuity equation \eqref{1.1}$_1$ and \eqref{123}, we can estimate $\|\nabla\rho_t\|_{L^2}$ as
\begin{equation*}
\begin{split}
\|\nabla\rho_t\|_{L^2}
& \leq C\|\nabla\mathbf{u}\|_{L^\infty}\|\nabla\rho\|_{L^2}
+C\|\rho\|_{L^\infty}\|\nabla^2\mathbf{u}\|_{L^2}
+C\|\mathbf{u}\|_{L^\infty}\|\nabla^2\rho\|_{L^2} \\
& \leq C\left(1+\|\nabla\mathbf{u}\|_{L^\infty}\right),
\end{split}
\end{equation*}
which together with \eqref{5.3} leads to
\begin{align}\label{5.4}
& \frac{d}{dt}\int\left(\rho|\nabla\mathbf{u}_{t}|^2
+|\nabla\mathbf{H}_{t}|^2\right)dx
+C_8\int\left(|\nabla^2\mathbf{u}_{t}|^2+
|\nabla^2\mathbf{H}_{t}|^2\right)dx\nonumber\\
& \leq CC_{3}(t)\left(1+\|\sqrt{\rho}\nabla\mathbf{u}_{t}\|_{L^2}^2
+\|\nabla\mathbf{H}_{t}\|_{L^2}^2\right).
\end{align}
Due to $\int_{0}^{T}C_{3}(t)dt<\infty$, we thus derive the desired \eqref{zx} from the Gronwall inequality.

We now estimate $L^{\infty}(0,T;H^3)$-bound of the strong solution $(\rho,\mathbf{u},\mathbf{H})$.
Similarly to \eqref{z10}, by standard elliptic estimates, we can deduce that
\begin{equation}\label{l3}
\|\mathbf{u}\|_{H^4}^2\leq C_{4}(t)\left(1+\|\nabla^3\rho\|_{L^2}^2\right),
\end{equation}
where $C_4(t)>0$ and $\int_{0}^{T}C_4(t)dt<0$. Then, taking the operator $\frac{\partial^3}{\partial x_i\partial x_j\partial x_k}$ on both side of \eqref{1.1}$_1$, and then multiplying the resulting equations by $\frac{\partial^3\rho}{\partial x_i\partial x_j\partial x_k}$, similarly to \eqref{10}, we can infer that
\begin{align}
\frac{d}{dt}\int|\nabla^3\rho|^2dx
 \leq C_5(t)\left(1+\|\nabla^3\rho\|_{L^2}^2\right),
\end{align}
where $C_5(t)>0$ and $\int_{0}^{T}C_5(t)dt\leq C$.
Hence we conclude by Gronwall's inequality that
\begin{equation}\label{l4}
\sup_{0\leq t\leq T}\|\nabla^3\rho\|_{L^2}\leq C,
\end{equation}
which together with \eqref{123} leads to
\begin{equation}\label{l5}
\sup_{0\leq t\leq T}\|\rho\|_{H^3}\leq C.
\end{equation}

We now estimate $\|\mathbf{u}\|_{H^3}$. It follows from \eqref{z3} and Lemma \ref{lem23} that
\begin{align}\label{l2}
\|\nabla^3\mathbf{u}\|_{L^2}
& \leq C\|\rho\partial_{t}\mathbf{u}
+\rho\mathbf{u}\cdot\nabla\mathbf{u}
+\nabla\rho^\gamma
-\nabla\nu\cdot\nabla\mathbf{u}
-\nabla\mathbf{u}\cdot\nabla\nu
-\divv\mathbf{u}\nabla\lambda
-(\nabla\times\mathbf{H})\times\mathbf{H}\|_{H^1}
\nonumber\\
& \leq C\|\rho\partial_{t}\mathbf{u}\|_{H^1}
+C\|\rho|\mathbf{u}||\nabla\mathbf{u}|\|_{H^1}
+C\|\nabla\rho^{\gamma}\|_{H^1} \nonumber\\
& \quad
+C\||\nabla\nu||\nabla\mathbf{u}|\|_{H^1}
+C\||\nabla\lambda||\nabla\mathbf{u}|\|_{H^1}
+C\||\mathbf{H}||\nabla\mathbf{H}|\|_{H^1}
 \nonumber\\
& \triangleq\sum_{j=1}^{6}L_j.
\end{align}
By \eqref{3.2} and \eqref{4.6}, we see that
\begin{equation*}
L_1\leq C+C\|\nabla\mathbf{u}_t\|_{L^2}.
\end{equation*}
From \eqref{123} and H{\"o}lder's inequality, it is not hard to check
\begin{equation*}
L_2+L_3+L_6\leq C.
\end{equation*}
The remainder $L_4$ and $L_5$ can be bounded as follows
\begin{equation*}
\begin{split}
L_4+L_5 & \leq C+C\|\nabla(|\nabla\nu||\nabla\mathbf{u}|)\|_{L^2} \\
& \leq C+C\||\nabla^2\rho||\nabla\mathbf{u}|\|_{L^2}
+C\||\nabla\rho||\nabla^2\mathbf{u}|\|_{L^2} \\
& \leq C+C\|\nabla^2\rho\|_{L^6}\|\nabla\mathbf{u}\|_{L^3}
+C\|\nabla\rho\|_{L^6}\|\nabla^2\mathbf{u}\|_{L^3}
 \\
& \leq C+C\|\nabla^3\rho\|_{L^2}
+C\|\nabla^2\mathbf{u}\|_{L^2}^{\frac{1}{2}}
\|\nabla^2\mathbf{u}\|_{L^6}^{\frac{1}{2}}
 \\
& \leq C+\frac{1}{2}\|\nabla^3\mathbf{u}\|_{L^2}.
\end{split}
\end{equation*}
Substituting the above estimates on $L_j$ into \eqref{l2}, and using \eqref{3.2} and \eqref{zx}, we have
\begin{equation}\label{l6}
\sup_{0\leq t\leq T}\|\nabla^3\mathbf{u}\|_{L^2}
\leq C+C\sup_{0\leq t\leq T}\|\nabla\mathbf{u}_t\|_{L^2}\leq C,
\end{equation}
which combined with \eqref{3.21} leads to
\begin{equation}\label{l7}
\sup_{0\leq t\leq T}\|\mathbf{u}\|_{H^3}\leq C.
\end{equation}

Finally, we bound $\|\mathbf{H}\|_{H^3}$. By standard elliptic estimates, it follows from \eqref{1.1}$_3$, \eqref{123}, \eqref{l7}, and \eqref{zx} that
\begin{equation}\label{l8}
\sup_{0\leq t\leq T}\|\mathbf{H}\|_{H^3}^2\leq C+C\sup_{0\leq t\leq T}\|\nabla\mathbf{H}_t\|_{L^2}^2\leq C.
\end{equation}

Therefore, we obtain the desired \eqref{l} from \eqref{l5}, \eqref{l7}, and \eqref{l8}. This finishes the proof of Lemma \ref{lem36}.
\hfill $\Box$

With Lemmas \ref{lem31}--\ref{lem36} at hand, we are now in a position to prove Theorem \ref{thm1.1}.

\textbf{Proof of Theorem \ref{thm1.1}.}
We argue by contradiction. Suppose that \eqref{B} were false, that is, \eqref{3.1} holds. Note that the general constant $C$ in Lemmas \ref{lem31}--\ref{lem36} is independent of $t<T^{*}$, that is, all the a priori estimates obtained in Lemmas \ref{lem31}--\ref{lem36} are uniformly bounded for any $t<T^{*}$. Hence, the function
\begin{equation*}
(\rho,\mathbf{u},\mathbf{H})(T^{*},x)
\triangleq\lim_{t\rightarrow T^{*}}(\rho,\mathbf{u},\mathbf{H})(t,x)
\end{equation*}
satisfy the initial condition \eqref{A} at $t=T^{*}$.
Therefore, taking $(\rho,\mathbf{u},\mathbf{H})(T^{*},x)$ as the initial data, one can extend the local strong solution beyond $T^{*}$, which contradicts the maximality of $T^{*}$. Thus we finish the proof of Theorem \ref{thm1.1}.
\hfill $\Box$

\section*{Appendix}

\textbf{Proof of \eqref{9.1}.} It follows from the mass equation \eqref{1.1}$_1$ that $\|\nabla\rho\|_{L^p}$ satisfies for any $2< p\leq6$,
\begin{align}\label{s}
\frac{d}{dt}\|\nabla\rho\|_{L^p}
\leq C\left(1+\|\nabla\mathbf{u}\|_{L^\infty}^2\right)\|\nabla\rho\|_{L^p}
+C\|\nabla^2\mathbf{u}\|_{L^p}.
\end{align}
To control the last term in \eqref{s},
we obtain from Lemma \ref{lem23}, \eqref{x2}, \eqref{x3}, Lemma \ref{lem22}, \eqref{3.2}, \eqref{3.7}, and \eqref{3.5} that
\begin{align}\label{s2}
\|\nabla^2\mathbf{u}\|_{L^p}
& \leq C\|\rho\|_{L^\infty}\|\nabla\mathbf{u}_{t}\|_{L^2}+
C\|\rho\|_{L^\infty}\|\nabla\mathbf{u}\|_{L^\infty}\|\mathbf{u}\|_{L^p}
+C\|\rho\|_{L^\infty}^{\gamma-1}\|\nabla\rho\|_{L^p}\nonumber\\
& \quad +C\|(\nu',\lambda')\|_{L^\infty}\|\nabla\mathbf{u}\|_{L^\infty}
\|\nabla\rho\|_{L^p}
+C\||\mathbf{H}||\nabla\mathbf{H}|\|_{L^p}
\nonumber\\
& \leq C\left(1+\|\nabla\mathbf{u}_{t}\|_{L^2}^2
+\|\nabla\mathbf{u}\|_{L^\infty}^2\right)\|\nabla\rho\|_{L^p}
+C\|\sqrt{|\mathbf{H}|}|\nabla\mathbf{H}|\|_{L^2}^{\frac{8-p}{3p}}
\|\sqrt{|\mathbf{H}|}\|_{L^8}^{\frac{4p-8}{3p}}
\nonumber\\
& \leq C\left(1+\|\nabla\mathbf{u}_{t}\|_{L^2}^2
+\|\nabla\mathbf{u}\|_{L^\infty}^2+\|\nabla^2\mathbf{H}\|_{L^2}^2
+\|\sqrt{|\mathbf{H}|}|\nabla\mathbf{H}|\|_{L^2}^2\right)
\left(1+\|\nabla\rho\|_{L^p}\right).
\end{align}
Substituting \eqref{s2} into \eqref{s}, we obtain that
\begin{align}\label{s3}
\frac{d}{dt}\|\nabla\rho\|_{L^p}
\leq CB(t)\left(1+\|\nabla\rho\|_{L^p}\right),
\end{align}
where
\begin{equation*}
B(t)\triangleq 1+\|\nabla\mathbf{u}_{t}\|_{L^2}^2
+\|\nabla\mathbf{u}\|_{L^\infty}^2+\|\nabla^2\mathbf{H}\|_{L^2}^2
+\|\sqrt{|\mathbf{H}|}|\nabla\mathbf{H}|\|_{L^2}^2 \in L^{1}(0,T)
\end{equation*}
due to \eqref{3.1}, \eqref{4.6}, \eqref{3.7}, and \eqref{3.5}. Hence we derive from the Gronwall inequality that
\begin{equation*}
\sup_{0\leq t\leq T}\|\nabla\rho\|_{L^p}\leq C.
\end{equation*}
This finishes the proof of \eqref{9.1}.   \hfill $\Box$

\end{document}